\documentclass{article}[14pt]
\usepackage{amsfonts}
\usepackage{amsmath}
\usepackage{amscd}
\newtheorem{thm}{Theorem}[section]

\newtheorem{defn}{Definition}[section]

 \numberwithin{equation}{section}
\newcommand{\RNum}[1]{\uppercase\expandafter{\romannumeral #1\relax}}
\begin{document}
\title{ Toeplitz Quantization  and Convexity }
\author{O. El Hadrami, Mohamed Lemine,  }
\date{Department of Mathematics \\  King Khalid Universtiy, Abha,
September 14, 2017} \maketitle
\begin{abstract}
Let $T^m_f $ be the Toeplitz quantization of a real $ C^{\infty}$ function defined on the sphere $ \mathbb{CP}(1)$.
$T^m_f $ is therefore a Hermitian matrix with spectrum $\lambda^m= (\lambda_0^m,\ldots,\lambda_m^m)$.
Schur's theorem  says  that the diagonal of a Hermitian matrix  $A$ that has the same spectrum of $ T^m_f $ lies inside a finite dimensional  convex set
whose extreme points  are $\{( \lambda_{\sigma(0)}^m,\ldots,\lambda_{\sigma(m)}^m)\}$, where
$\sigma$ is any permutation of $(m+1)$
elements. In this paper, we  prove that these convex sets "converge" to a huge convex set in $L^2([0,1])$  whose extreme points are $ f^*\circ \phi$, where $ f^*$ is the decreasing rearrangement of $ f$ and $ \phi $ ranges over the set of  measure preserving transformations of the unit interval $ [0,1]$.
 \end{abstract}
\section{Introduction and background}
In their papers \cite{BFR,BHFR}, the authors have described  similarities between the infinite  and finite  dimensional Lie groups. They have strengthen  the idea that  the set  $SDiff(\mathbb{CP}(1))$ of area preserving diffeomorphisms of the Riemann sphere
is an infinite dimensional analog of $SU(n)$ by
 proving an infinite version  of the $ SU(n)$ Schur-Horn convexity theorem .\\
In the present paper, we want to show that the convex sets in the two versions (finite and infinite) of  Schur-Horn convexity theorem are related.
 In order to do that, we will show first that each permutation of $(m+1)$ letters determine a measure preserving transformation of the interval $[0,1]$.
 Secondly we use the fact that the eigenvalues of the Toeplitz quantization of $ f$ determine the decreasing rearrangement of $f$.\\
   But since the two convex sets are defined by inequalities arising from   the theory of majorization in $ \mathbb{R}^{m+1}$ developed in \cite{Ra,Rak}  and its generalization  to $L^1([0,1])$ by J.Ryff \cite{Ry1,Ry2,Ry3},   we start by summarizing  briefly here the main points.\\
Majorization is a partial ordering in  $\mathbb{R}^{m+1}$ defined as it follows:\\
For $ x \in \mathbb{R}^{m+1},$ let $ x^* $ denote the vector obtained by rearranging the components of $ x$ in non-increasing order. We say that $ x$ majorizes $ y,$ written $ y \prec x, $ if
$$ y_0^* + y_1^* + \cdots + y_k^* \leq x_0^* + x_1^* + \cdots + x_k^* , \quad 0 \leq k \leq m-1$$
$$ \sum_{k=0}^{k=m}y_k^* =  \sum_{k=0}^{k=m}x_k^*$$
Now we can state Schur' theorem for hermitian matrices.
\begin{thm}\label{t:sch}
 Let $  \lambda^m = ( \lambda_0,\cdots,\lambda_m) \in \mathbb{R}^{m+1}$ be the eigenvalues of a hermitian matrix  $A$ . Let $ diag(A)= ( a_{00}, \cdots, a_{mm})$ be the diagonal of $A$
  then $ \lambda^m $ majorizes $ diag(A)$.
  \end{thm}
  i.e
 \begin{equation}\label{eq:maj}
 diag(A) \prec \lambda^m
\end{equation}
Before going on, let us first fix some notations:
\begin{enumerate}
  \item Let $ x = (x_0, \cdots, x_m) \in \mathbb{R}^{m+1}.$ \\
  $\sum_m x$  is the orbit of $ x$ under the symmetric group of $ (m+1) $ letters, i.e the collection of points $(x_{\sigma(0)},\cdots,x_{\sigma(m)}),$ where $ \sigma $ ranges over all $ ( m+1)!$ permutations.
    \item For $ C \subset \mathbf{E}$ where $\mathbf{ E}$ is a vector space over $ \mathbb{R}, co(\mathbf{C})$ denote the convex hull of $\mathbf{C}.$
\end{enumerate}
Majorization and convexity are closely related as it is shown by the  following theorem.
 \begin{thm}\label{t:rado} Let $ x \in \mathbb{R}^{m+1}$
\begin{itemize}
  \item ( Rado's theorem)
  $$  \{ y \in \mathbb{R}^{m+1}, y \prec x \}  = co ( \sum\nolimits_m x)$$
  \item $y \prec x $  if and only if $ \sum_{i=0}^{i = m} f(y_i) \leq  \sum_{i=0}^{i = m} f(x_i)$ for any convex function whose domain contains all the numbers $ x_i, y_i, 0 \leq i \leq m$.
\end{itemize}
\end{thm}
 Schur-Horn' theorem  can  be therefore  restated in the following terms :
\begin{thm}[Schur-Horn's theorem] \label{t:SC}
Let $ \lambda^m = ( \lambda_0,\cdots,\lambda_m) \in \mathbb{R}^{m+1}.$\\
 Let $ \mathcal{H}_{\lambda^m }$ be the set of all hermitian matrices whose spectrum is
$ \lambda^m.$ \\
Let $ p^m :\mathcal{H}_{\lambda^m } \longrightarrow \mathbb{R}^{m+1}$ be the map that picks the diagonal of a matrix.\\
 Then the image of the map $ p^m $ is the convex set $ co( \sum\nolimits_m{\lambda^m}).$
\end{thm}
The convex set $ co( \sum\nolimits_m{\lambda^m})$ plays a very important role in symplectic geometry: It is the image of a moment map \cite{A,GS}.\\\\
The concepts of majorization is also extended to integrable functions in the following sense.\\
   Let $(X,\mu) $ be a finite measurable space. For $ f $ measurable function on $ X$, the distribution function of $ f$ is the function $ F_f$ defined by
$$ F_f(t) = \mu( \{ \omega, f(\omega) < t\}) .$$
Let $ d_f(t) = \mu(X) - F_f(t)$.
\begin{defn}
The decreasing rearrangement of $ f$  is the function $ f^*$ defined by:
$$ f^*(s) = \inf \{ t > 0, d_f(t) < s\}.$$
\end{defn}
    For $f, g \in L^1((X,\mu))$, let $ f^*$ and $g^*$ be their  decreasing rearrangement \\
   We say that $  f $ majorizes $ g $ ( written $ g \prec f $ ) if
   \begin{align*}
   \int_0^s g^*(z)\,dz & \leq  \int_0^s f^*(z)\,dz, \quad 0 \leq s  < 1 \\
    \int_0^1 g^*(z)\,dz & = \int_0^1 f^*(z)\,dz
   \end{align*}
   To stay conform with the notation of \cite{BHFR}, let $X = \mathbb{CP}(1)$  be the Riemann sphere , $\mu $ the measure defined by the the Fubini-Study symplectic form  which in the local coordinate $ [1,w]$ is given by
$$ \Omega = \frac{i}{(1 + w\overline{w})^2}\,dw\wedge d\overline{w}.$$
   Set $ w = r\exp(i \theta)$ and introduce  the real variable $ z\in [0,1[$ by $z = r^2 \backslash(1+r^2)$. The symplectic form $ \Omega $ becomes
   $ \Omega = 2dz \wedge d \theta$.
The infinite version of Schur theorem is
\begin{thm}[Schur-type Theorems,\cite{BFR}]\label{t:Schur}
    Let $ L^2(\mathbb{CP}(1)) $ be the set of square integrable functions on the sphere.
        Let $ P: L^2(\mathbb{CP}(1)) \longrightarrow L^2[0,1] $ be the projection $\displaystyle{ P(f) = \frac{1}{2 \pi}\int_0^{2 \pi} f(z,\theta)d\theta}$.\\\\
    Then $ f^* $ majorizes $ P( f)$.
    \end{thm}
    i.e
$$ P(f) \prec f^*$$
We have also   the equivalent of Rado's theorem in $  \in L^2([0,1])$.
\begin{thm} Let $ f \in L^2([0,1])$. The set $ \displaystyle{\Omega(f) = \{g \in L^2([0,1]),  g \prec f \}}$ is weakly compact and convex.
Its set of extreme points is \\
 $\displaystyle{ \{ f^* \circ \phi\; |\; \phi \; \mbox{ is a measure preserving transformation of }\; [0,1] \}}$.
\end{thm}
The paper is organized in three sections: In $\S2$, we have reviewed  the topology of the set of measure preserving transformations of the unit interval$[0,1]$
 and showed that  the groups, $ \sum_m$ of permutations of $ (m+1)$ letters, can be identified with a dense subset of the
  set of all invertible measure preserving transformations of $ [0,1]$ endowed with strong operator topology.\\
    In $\S3$ we use Toeplitz quantization to show that  the weak closure of the topological lim sup  of the sets $ co\left(\sum_m\cdot \lambda^m\right)$ is the set \\
    $\Omega(f^*) = co \left( \{ f^* \circ \phi, \phi \;\mbox{ measure preserving transformation of } \; [0,1]\}\right)$.
           \section{measure preserving transformations of [0,1]}
  The Lebesgue measure on the unit interval $ I = [0,1]$ will be always denoted by $ \lvert\, \cdot \, \rvert$.\\
 A map $ \phi$ from $ [0,1]$ to itself is a measure preserving transformation if
$$  \lvert \phi^{-1} (A) \rvert =  \lvert  A \rvert,\qquad \mbox{for Borel set A}  $$
The set of all (non necessary invertible ) measure preserving transformation of the unit interval will be denoted \emph{Smeas(I)}.
 The invertible ones will be denoted by Imeas(I).\\
    Each   $ S  \in Smeas(I) $ determine  a bounded linear operator $ P_S $ on  $ L^2([0,1])$ by
   $ P_S (f) = f \circ S$. In this way, $ Smeas(I)$ can be identified to a subset of the set of bounded linear operators of $ L^2[0,1])$  and the strong operator topology induces a topology on $Smeas(I)$.\\
   Evidently, a sequence  $ S_n $ converges to $ S$ in the strong operator topology if  for every function $f,\,  f \circ S_n $ converges   to  $ f \circ S $ in $ L^2([0,1])$.\\
  \\\\
To state our first main result, we need to define dyadic permutations.\\
    Let  $I_k^m =  \left[  k \backslash (m+1) , (k +1)\backslash(m+1)\right), \, k =0, 1, \cdots, m;\, m = 0, 1, \cdots $.\\
     Let $\sum_m$ be the group of permutations of $(m+1)$ letters. For $\sigma \in \sum_m$,
     $ \widehat{\sigma}$ is the invertible measure preserving transformation that sends the interval $ I^m_k$ to the interval $I^m_{\sigma(k)}$ by ordinary translation.\\
    We call $ \widehat{\sigma}$ a  permutation  of rank $m$. In this way, we can identified the group of permutations $ \sum_m$ with a subgroup of Imeas(I).\\
      Halmos (\cite{Ha})  shows that the set of all permutations of different rank is dense in Imeas(I) for the strong operator topology.\\
       Also  Brown, in \cite{Br} has proved that $Smeas(I)$ is the closure of $ Imeas(I)$ for the strong operator topology.\\
 We can summarize these results in the following theorem.
 \begin{thm}\label{t:rsg} Let $ \sum_{m}$ be the symmetric group of $(m+1)$ letters. There exists a one to one  group homomorphism
$ \Psi^m :\sum_{m} \longrightarrow \emph{Imeas(I)}$ that sends $ \sigma$ to $ \widehat{\sigma}$ and if we identify $ \sum_m$ with $\Psi^m(\sum_m)$, then $ \bigcup_m\sum_m $ is a dense set in Smeas(I) for the strong operator topology.
  \end{thm}
  \section{Toeplitz Quantization and Convexity}
  \subsection{Toeplitz Quantization }  Consider the Riemann sphere $ \mathbb{CP}(1)$ with
   the Fubini-Study symplectic form  in the local coordinate $ [1,w]$
$$ \Omega = \frac{i}{(1 + w\overline{w})^2}\,dw\wedge d\overline{w}$$
and the standard hyperplane bundle $\mathbf{L}$.
    The tensor power  $ \mathbf{L}^{\otimes m} $  has $ (m +1)$ linearly independent sections which in the local coordinate $w$ are just $ 1, w,...,w^m$.  The
    bundle   $ \mathbf{L}^{\otimes m} $ comes equipped with the Hermitian metric
    $$ \langle s_1,s_2\rangle(w) = \dfrac{1}{(1 + \lvert w \rvert^2)^m}s_1(w)\overline{s_2(w)}.$$
    Now let $\Gamma^m_2$ be the space of square-integrable sections of $ \mathbf{L}^{\otimes m} $ and
   $ \Gamma_{hol}^m $  the space  of holomorphic sections ( the span of  $ 1, w,...,w^m$).The orthogonal projection $\Gamma^m_2 \longrightarrow \Gamma_{hol}^m $
   is denoted by $ P^m$.The Toeplitz quantization of $ f$  is the map $ T_f^m\,: \Gamma^m_{hol} \rightarrow \Gamma^m_{hol}$ defined by
$$T_f^m = P^m \circ \mathbf{M}_f\circ P^m $$ where $\mathbf{M}_f $ is multiplication by $f$.
We refer the interested reader to \cite{BMS} for a detailed exposition on Toeplitz quantization.\\
The crucial result is the following theorem
\begin{thm}[Distribution of the Eigenvalues of Toeplitz Quantization ]\label{t:Had}
Let  $  \lambda^m = \left(\lambda_0^m , \lambda_1^m , \cdots ,\lambda_m^m  \right)$ be the eigenvalues of $ T_f^m $ arranged in non-increasing order and let  $\Lambda^m (s) $  be the real step function defined on the interval $[0,1[$ by
\begin{equation}\label{e:lamda}
   \Lambda^m \left(\left [ \frac{k}{m+1}, \frac{k+1}{m+1}\right[ \right) = \lambda_k^m, \quad  0 \leq k \leq m .
   \end{equation}
  Then  the sequence  $ \left(\Lambda^m(s)\right)_m$  converges  point-wise
almost everywhere  to the decreasing rearrangement $ f^*(s) $ of the function $f$.
\end{thm}
 The proof of this Theorem is based  on the the following theorem
\begin{thm}[Szeg\"{o}-type Theorem,\cite{Gu} p: 248]\label{t:Gu} Given a smooth real valued function $f$ on $\mathbb{CP}(1)$, let $T^m_f $ be Toeplitz quantization of $f$ and  let $ \mu^m$ be its  spectral measure . Then
   $\dfrac{\mu^m}{m+1}$ tends weakly to a limiting measure as $ m$ tends to infinity, this limiting measure being
   $$   \mu(\phi) = \frac{1}{2\pi}\int_{\mathbb{CP}(1)} \phi(f(x)) d\Omega, \quad \mbox{for}\, \phi \in \mathbf{C}(\mathbb{R}).$$
        \end{thm}
i.e if $ \left(\lambda_0^m , \lambda_1^m , \cdots ,\lambda_m^m  \right)$ are the eigenvalues of $ T_f^m $ ordered in non-increasing order, then
\begin{equation}\label{e:szego1}
\lim_{m\longrightarrow +\infty} \sum_{k=0}^{k=m}\dfrac{\phi(\lambda_k)}{m+1} = \frac{1}{2\pi}\int_{\mathbb{CP}(1)} \phi(f(x)) d\Omega.
\end{equation}
If we use the step function $ \Lambda^m$ defined be (\ref{e:lamda}) then (\ref{e:szego1}) can be written
\begin{equation}\label{e:sz}
 \lim_{m\longrightarrow +\infty} \int_0^1 \phi(\Lambda^m) = \frac{1}{2\pi}\int_{\mathbb{CP}(1)} \phi(f(x)) d\Omega .
 \end{equation}
  But since $ f$ and $f^*$ are equi-measurable, we have
 $$  \frac{1}{2 \pi}\int_{\mathbb{CP}(1)} (\phi \circ f) d\Omega = \int_0^1 \phi \circ f^*(t) dt $$
 and (\ref{e:sz}) becomes
\begin{equation}\label{eq:spc}
   \lim \limits_{ m \mapsto +\infty}\int_0^1 \phi(\Lambda^m)(t)\, dt  = \int_0^1 \phi (
f^*)(t) \,dt.
\end{equation}
where $\Lambda^m$ is defined by (\ref{e:lamda}).\\
Relation ( \ref {eq:spc} ) is equivalent to:  the sequence of step functions $ \Lambda^m$
converges in distribution to the real function $f^*.$ (One can see ~\cite{Du} for more details on convergence in distribution.)\\
 In general  convergence in distribution does not imply convergence point-wise. Nevertheless,  there exists another sequence $g_n$ with the same distribution as $f_n$ that
converges point-wise to a function $g$, that has the same distribution of $ f$. That is the content of Skorokhod's Theorem.
\newpage
\begin{thm}[Skorokhod]\label{t:SK}
   Let $ (X,\Sigma,\mu)$ be a finite measure space,and $f,\,f_n:\,X \rightarrow \mathbb{R}$
   be a sequence of measurable functions such that $f_n$ converge in distribution to $f$. Then on the
   Lebesgue  measure space $ ( I, \mathcal{B},\left\lvert \cdot \right\rvert)$, where $ I = (0, \mu(X), \mathcal{B}$ is the Borel $\sigma$-algebra of $I$, and $ \left \lvert \cdot \right\rvert$ is the Lebesgue measure, there exists measurable functions  $ g_n,g: I \rightarrow \mathbb{R}$ such that $ g_n(t) \rightarrow g(t) \;a.e [\left\lvert \cdot \right\rvert]$, and
   \begin{align*}
    \mu ( \{\omega: f_n(\omega) < x \}) &= \left \lvert  \{ t: g_n(t) < x\} \right \rvert \\
    \mu ( \{\omega: f(\omega) < x \}) &= \left \lvert \{ t: g(t) < x\} \right \rvert
    \end{align*}
   $ x \in \mathbb{R}, n \geq 1$.
\end{thm}
Let $F_n $ and $ F$ be the distribution functions of $f_n$ and $f$.
We can take $g_n$ and $ g$  to be just the generalized inverse of $F_n$ and $F$:
$$g_n(t) = \inf\{x: F_n(x) > t\},\quad g(t) = \inf\{x: F(x) > t \}, \; 0 < t < 1 .$$
It easily seen then that
$$ g_n(1 -t) = g_n^*(t) = f_n^*(t), \quad g(1 -t) = g^*(t) = f^*(t).  $$
(  \cite{Ra} [p141-144] .)\\
  Applying  Skorokhod's Theorem to the the sequence $ \Lambda^m $, we deduce that
the generalized inverses of the  functions  $ \Lambda^m$ converges point-wise almost everywhere
to the generalized inverse of the function $f$.\\
Consequently , the sequence of the decreasing rearrangements of $\Lambda^m $ converges to the decreasing rearrangement  $ f^* $ of $ f$.
$$ \lim_{m \longmapsto +\infty}\left(\Lambda^m \right)^*(s) = f^*(s).$$
But since $\Lambda^m(s)$ is a decreasing function,  $\left(\Lambda^m\right)^*(s) = \Lambda^m(s)$, and therefore we have
$$  \lim_{m \longmapsto  +\infty}\Lambda^m(s) = f^*(s) $$ and that
 ends  the proof of the theorem.

\subsection{Convexity}
  Let    $ p^m: \mathbb{R}^{m+1} \rightarrow L^2([0,1])$  be the map that associates to the point $(a_0,a_1,\cdots,a_m)$ the step function
   $ p^m(a_0,a_1,\cdots,a_m) = \sum_{k=0}^{k=m}a_k\chi_{I^m_k}$, where $I^m_k$ is the interval $\left [\dfrac{k}{m+1},\dfrac{k+1}{m+1}\right [$
   and $  \chi_{I^m_k}$ is the characteristic function of $I^m_k$. Set
    $ E_m = p^m( co(\sum_m \cdot \lambda^m))$.\\
   Now we are ready to state the main result about convexity.
\begin{thm}Let $ f \in C^{\infty}(\mathbb{CP}(1))$. Let $ f^* $ be the decreasing rearrangement of $ f$.  Let
  $T_f^m $ be the Toeplitz quantization of $ f$ and  let  $\lambda^m   = (\lambda^m_0,\cdots,\lambda^m_m)$ be
   the eigenvalues of $ T_f^m$. Let $\Omega(f^*) = \{ g \in L^2([0,1], g \prec f^*\}$.
   Then the closed convex hull of the set $ Smeas(I) \cdot f^* = \{f^* \circ \phi \;| \phi \; \in Smeas(I)\}$ is the weak closure of the topological lim sup  of  the convex sets  $ E_m = p^m( co(\sum_m \cdot \lambda^m))$.
   \end{thm}
We recall the definition of the closed limit.
\begin{defn}[ \cite{AB}, p. 114]
Let $ \{ E_m \}$ be a sequence of subsets of a topological space X.
Then,  a point x in X belongs to \textbf{the topological lim sup} of $ E_m$, denoted  $ Ls E_m $, if for every
neighborhood V of x there are infinitely many m with\\ $ V \cap E_m  \neq \emptyset. $
 \end{defn}
 Clearly,  $Ls E_m$ is a  closed sets  and moreover
$$Ls E_m  = \bigcap_{m=1}^{m= \infty} \overline{\bigcup_{k=m}^{m= \infty} E_k}.$$
In our case, $ X = L^2([0,1])$  and $ E_m = p^m( co(\sum_m \cdot \lambda^m))$ .\\
 Since  $ L^2([0,1], \lVert \cdot\rVert_2)$ is a normed vector space, every point of $ L^2([0,1] $ has a countable basis of neighborhood and    $ LsE_m$ can be characterized in terms of sequences:\\
    $ g \in LsE_m $ if and only if there exists a subsequence $ (g_{m_k})_k $ such that $ g_{m_k} \in E_{m_k}$ and $ g_{m_k}$ converges to $ g \in L^2([0,1])$.
        \\
  Now we claim
   \begin{description}
              \item[(A)]  $Smeas(I) \cdot f^* \subset LsE_m.$
                             \item[(B)] $ \Omega(f^*)  \subset \overline{LsE_m}^{weak} \subset \Omega(f^*)$.
                          \end{description}

\textbf{Proof of( A):}\\
  Let $ f^* \circ \phi \in Smeas(I) \cdot f^* $.\\
  The set of dyadic permutations is a countable set and  from theorem (\ref{t:rsg}) is dense in $Smeas(I)$  for the strong operator topology; therefore  there exists a sequence  $(\sigma_n)_n$  of  permutations that converges to  $\phi$ , i.e

$$ \forall  g \in L^2([0,1]), \forall \epsilon >0, \exists N_0, \forall n \geq N_0, \| g \circ \sigma_n - g \circ \phi \| _2 < \dfrac{\epsilon}{2}. $$
  In particular for $ g =\Lambda^m$, we have,
\begin{equation}\label{e:con1}
 \forall  m , \forall \epsilon >0, \exists N_0,  \forall n \geq N_0, \lVert \Lambda^m\circ \sigma_n - \Lambda^m \circ \phi \rVert _2 < \dfrac{\epsilon}{2} .
 \end{equation}
 Now since  the sequence $\Lambda^m $ converges almost everywhere  to $ f^* $ and $\forall  m, \forall x \in [0,1], \left| \Lambda^m (x)\right| \leq  \lVert f^* \rVert_{\infty},, $ it converges also to $f^* $ in $L^2([0,1])$, ie
 \begin{equation}\label{e:con2}
 \forall \epsilon > 0, \exists M_0,  \forall m \geq M_0, \lVert \Lambda^m\circ \phi   - f^* \circ \phi \rVert _2 < \dfrac{\epsilon}{2}.
\end{equation}
We conclude then from (\ref{e:con1}) and (\ref{e:con2}) that\\
$\forall \epsilon >0, \exists M_0,\exists N_0, \forall n \geq N_0, \forall m \geq M_0$
\begin{equation}\label{e:con3}
             \| \Lambda^m \circ \sigma_n - f^* \circ \phi \| _2 \leq \| \Lambda^m \circ \sigma_n - \Lambda^m \circ \phi \| _2 + \| \Lambda^m \circ \phi  - f^* \circ \phi \| _2  \leq \epsilon.
\end{equation}
Now if $ \sigma_n $ is of order $ k_n$, and if  we let $ m = k_n$ in (\ref{e:con3}) we get

      $ \forall \epsilon >0, \exists M_0,\exists N_0 $ for every $ n \geq N_0 $ such that $ k_n \geq M_0$, we have
      $$\lVert \Lambda^{k_n} \circ \sigma_n - f^* \circ \phi \rVert _2  \leq \epsilon.$$
      Since  $\Lambda^{k_n} \circ \sigma_n \in E_{k_n} $, we conclude
    therefore that $ f^* \circ \phi \in  LsE_m$ .\\
    \textbf{Proof of (B)}:\\
    It is shown in [~\cite{BFR},p 523] that $ \Omega(f^*) = \{ g \in L^2([0,1]), g \prec f^* \} $ is weakly compact and convex. Its extreme points are precisely  the elements $ f^* \circ \phi , \phi \in Smeas(I)$. It  is also indicated in  [~\cite{Ry3},~ p1030] that the set of extreme points is dense in $ \Omega(f^*)$ for the weak topology. But we have just seen in part (A) of our claim  that $ Smeas(I)\cdot f^* \subset LsE_m$ . Taking the closure in the weak topology we get
       $$ \Omega(f^*) = \overline{Smeas(I)\cdot f^*}^{weak} \subset \overline{LsE_m}^{weak}.$$
     It remains to show that $  \overline{LsE_m}^{weak} \subset \Omega(f^*)$.\\
      Let $ g \in LsE_m$. Then by definition, there exists a subsequence $ (g_{m_k})_k$ such that
     $ g_{m_k} \in E_{m_k}, $ and $(g_{m_k})_k$ converges in  $L^2([0,1])$ to $ g.$\\
     Now by Rado's theorem (\ref{t:rado}) , we have
   $$ g_{m_k} \in E_{m_k} \Longleftrightarrow  g_{m_k} \prec \Lambda_{m_k} $$
   But the sequence  $\left(g_{m_k}\right)_k$  being convergent in  $L^2([0,1])$), we can extract a subsequent $\left(g_{m_{k_l}}\right)_l$ that converges pointwise to $ g$ a.e.
   We have then
   $$  g_{m_{k_l}} \prec  \Lambda_{m_{k_l}}.$$
   And taking limit (simple convergence)  of both sides, we get
   $ g \prec f^* $  i.e  $ LsE_m \subset \Omega(f^*)$. Taking the weak closure of both sets
 we conclude therefore that
 $$ \Omega(f^*) =  \overline{LsE_m}^{weak}.$$
\newline
  Our next goal is to compare   the co-adjoint orbits of $SU(m+1)$ and the co-adjoint orbits of SDiff($\mathbb{CP}(1)$, the group of area preserving diffeomorphisms of the sphere.
 \section*{Acknowledgment}
     I would like thank  Professor H. Flaschka  for valuable discussions.
  
King Khaled University, Abha, Saudi Arabia.\\
E-mail address: mbouleryah@kku.edu.sa
\end{document}